# On graph coloring check-digit method

Kamil Kulesza  Zbigniew Kotulski

Institute of Fundamental Technological Research,
Polish Academy of Sciences,
Świętokrzyska 21, 00-049 Warsaw, Poland

September 11, 2002

## Abstract

We show a method how to convert any graph into the binary number and vice versa. We derive upper bound for maximum number of graphs, that, have fixed number of vertices and can be colored with $n$ colors ($n$ is any given number). Proof for the result is outlined. Next, graph coloring based check-digit scheme is proposed. We use quantitative result derived, to show, that feasibility of the proposed scheme increases with size of the number which digits are checked, and overall probability of digits errors.

## 1 Introduction

Development of modern digital computers resulted in growing popularity of discrete mathematics structures like graphs. Nowadays graphs are used in virtually every field of computer science. Problems of proper vertex coloring for an arbitrary graph, with minimal set of colors are known to be NP complete, see [1], [2]. The same is true for decisions problems about the graph coloring (e.g., on the graph's chromatic number). Coloring and the chromatic number are integral properties of any graph. Given the graph $G$, it is always possible to find its chromatic number and $n$-coloring. Hence, one may try to use them for own advantage.

We propose check-digit scheme based on the graph coloring. Due to the volume constrains, we present only basic ideas, just to illustrate the concepts. The outline of the paper is the following: in section 2 procedures to convert any graph into the number and number into the graph are stated. The section is summarized with short discussion about graph $n$-coloring. Next, we state theorem about upper bound on number of graphs, that have fixed number of vertices and can be colored with $n$ colors (section 3). In section 4 we describe our proposal for graph coloring based check-digit scheme. Section 5 provides conclusions and description of further research.

## 2 Preliminaries

At the beginning simple scheme that allows conversion between graphs and numbers is presented. While more sophisticated methods can be used (e.g. [5]), the one chosen well illustrates development of main results.

Let $G(V, E)$ be a graph, where $V$ is a set of vertices and $E$ is a set of edges, with $|E|$ edges and $|V|$ vertices; $v_i$ denotes $i$th vertex of the graph, $v_i \in V$. $K_n$ is a complete graph on $n$ vertices, $\chi(G)$ is chromatic number for the graph $G$.

In further considerations the graph $G$ with $m$ vertices is described by the square adjacency matrix $\mathbf{A} = [a_{ij}]$, $i, j = 1, 2, ..., m$. The elements of $\mathbf{A}$ satisfy:

- for $i \neq j$, $a_{ij} = 1$ if $v_i v_j \in E$ (vertices $v_i, v_j$ are connected by an edge) and $a_{ij} = 0$, otherwise;

- for $i = j$, $a_{ii} = \alpha$, where $\alpha \in Z_k$ is the number of color assigned to $v_i$. In $Z_k$, $k \geqslant \chi(G)$ denotes the number of colors that can be used to color vertices of $G$ (in other words, $k$ is the size of the color palette).

In the case that the graph coloring is not considered, $k = 1$, and all entries on the main diagonal of $\mathbf{A}$ are zero.

$\mathbf{A}$ is a symmetric matrix, hence having all the entries on the main diagonal and all the entries below main diagonal, one can describe whole matrix (and as the result graph $G$). Thus it can be written as the sequence $a_{21}a_{31}a_{32}a_{41}a_{42}a_{43} ... a_{m(m-1)}a_{11} \ a_{22}...a_{mm}$, where the first binary part ($a_{21}a_{31}$ $a_{32}a_{41}a_{42}a_{43} ... a_{m(m-1)}$) corresponds to all the entries below main diagonal, while second decimal one ($a_{11} \ a_{22}...a_{mm}$) to the main diagonal itself.

Coding the number as a graph can be done by converting number into binary form. Then the adjacency matrix $\mathbf{A}$ is encoded in the same way as described above. If the length $l$ of binary number does not yield integer solution to the equation $l = \dfrac{(m-1)\,m}{2}$, then lacking matrix entries can be added. Method for entries addition is further discussed in the section 4.

Finding graph minimal coloring or even decision on the chromatic number are complex and difficult to comprehend. For instance: there are numerous examples that single edge change in the graph structure can result in dramatic change in its coloring. Any significant result in this field would bring us closer to resolving fundamental P versus NP question. Current research concentrates on three fields:

1. Algorithms. Building better approximation solution algorithms, that have polynomial time characteristics and allow to color graph with $\chi(G) = n$ with the smallest $k$-color palette ($|V| \geqslant k \geqslant n$). Results available allow to obtain quite reasonable $k$ values (e.g., see [3]). It is important to note that mentioned above algorithms can work in the deterministic way.

2. Probabilistic methods. Present state of the art is given in [2].

3. Theory. Finding new coloring criterions and new cases, that are solvable in polynomial time (e.g., see [4]). To give reader a flavor few well known results, mentioned later in the text, are presented.

**Theorem.** *Every graph $G$ of $K_n$ configuration has $\chi(G) = n$.*

**Lemma.** *Every graph $G$ having a subgraph of $K_n$ configuration has $\chi(G) \geqslant n$.*

**The Brook's theorem** (1941). *If the graph $G$ is not an odd circuit or a complete graph, then $\chi(G) \leqslant d$, where $d$ is the maximum degree of a vertex of $G$.*

In further considerations, two observations are used:

1. Vertex coloring of the graph can be checked only when the structure of the graph is known.

2. Checking whether a graph is properly colored can be done in linear time, see [1].

## 3 Upper bound on number of graphs with fixed $|V|$ and $n$

Let's denote $\Gamma(|V|)$ as the number of all possible graphs having $|V|$ vertices. One can calculate that

$$\Gamma(|V|) = 2^{\frac{|V|(|V|-1)}{2}}. \tag{1}$$

In other words, $\Gamma(|V|)$ represents number of binary sequences such that every bit corresponds to the possible edge in the $|V|$ vertices graph (1 denotes presence of the edge, 0 opposite case).

The graph coloring considered in this section is merely graph $n$-coloring ($n \geqslant \chi(G)$) not the minimal coloring of the graph ($n = \chi(G)$). Let's denote:

$x_j$ is the number of vertices colored by $j$th color ($0 < j \leqslant n$). Vertices are partitioned into $n$ sets, such that all vertices with the same color are in one set,

$V_j$ is the set of vertices colored by $j$th color ($0 < j \leqslant n$), $|V_j| = x_j$,

$P$ is a particular, non-degenerate partition of $|V|$ into $n$ pieces. $P = \{x_1, x_2, ..., x_j, ..., x_n\}$, $\sum_{j=1}^{n} x_j = |V|$ and $x_j > 0$ for all $j = 1, 2, ..., n$,

$\Gamma(|V|, n, P)$ is the number of graphs having $|V|$ vertices and particular partition $P$ ($|V|$ divided into $n$ pieces),

$\Gamma(|V|, n) = max_P \{\Gamma(|V|, n, P)\}$ where $max_P$ denotes the maximum value of $\Gamma(|V|, n, P)$ over of all possible partitions $P$ for given $|V|$ and $n$.

**Theorem (upper bound for $\Gamma(|V|, n)$),**

Let $|V| = n + y$ where $y \in N$. Then

$$\Gamma(|V|) \geqslant 2^y \Gamma(|V|, n). \tag{2}$$

**Outline of the proof:**

First, observe that

$$\Gamma(|V|, n, P) = 2^{(x_2+x_3+...+x_n)x_1 + (x_3+x_4+...+x_n)x_2 + ... + x_n x_{n-1}}. \tag{3}$$

This allows to carry out further work on $\log_2$ (base 2 logarithms)

There are two cases to be considered:

Case 1: $|V| = n$. Hence there is only one partition $P$ (all $x_j = 1$), substitution into (3) yields $\Gamma(|V|) = \Gamma(|V|, n)$ ($y = 0$, so $2^y = 1$).

Case 2: $|V| > n$. First properties of $\log_2(\Gamma(|V|, n, P))$ for all possible $n$ values are examined to find general formula.

Comparing corresponding terms from $\log_2(\Gamma(|V|))$ and $\log_2(\Gamma(|V|, n, P))$ one finds that

$$\log_2(\Gamma(|V|)) \geq y + \log_2(\Gamma(|V|, n, P))$$

for any $P$. When $|V| > 2n$ condition stated above holds with strict inequality for any $P$. Ultimately:

$$\Gamma(|V|) \geqslant 2^y \Gamma(|V|, n)$$

∎

**Remarks:**

1. Reasoning analogous as in the proof of 2 leads to results stating the relations between $\Gamma(|V|)$ and $\Gamma(|V|, n)$ for more specific conditions (say, it exists $j$ such that $x_j = k$, $k \geqslant 2$).

2. Technique used to prove (2) yields the number of the graphs that can be properly colored with $n$ colors. The result includes class of graphs $G'$ such that $\chi(G') < n$. This observation can be used to calculate the number of the graphs that have fixed chromatic number $\chi(G) = n$.

Obtained results are not crucial for check-digit scheme outline, hence they are omitted.

# 4 Check-digit scheme based on the vertex graph coloring

To simplify procedure's description, it is assumed that there is an efficient procedure for finding graph $G$ minimal $n$-coloring ($n = \chi(G)$). Issue of finding graph minimal coloring will be addressed at the end of this section.

Let us introduce the notation:

$D$ denotes number that is stored/sent, $D'$ denotes number that is read/received

$G(D), G(D')$ denotes graphs resulting from conversion of $D$, $D'$ into the graph form

$col(D), col(D')$ denotes minimal $n$-coloring for graphs $G(D), G(D')$, respectively

*PROCEDURE (checks whether $D = D'$)*:

Encoding

1. Convert number $D$ to graph $G(D)$.
2. Find the proper coloring of graph, $col(D)$.
3. Store/Send $\{D, col(D)\}$.

{Process that can modify $D$ (e.g., transmission) takes place}

Decoding
It is assumed that $col(D') = col(D)$[1] .
1. Read/Receive number $D'$ and convert it to graph $G(D')$
2. Is $col(D')$ valid $n$-coloring for $G(D')$?
Results: No→error found(terminate), Yes → continue
3. Is $G(D')$  $(n-1)$-colorable?
Result: Yes→error found(terminate), No → positive verification for $D' = D$
End of procedure.

***Analysis & Discussion***

Our check-digit scheme should be discussed in the context of suitable communication model. Nevertheless we provide few simple observations that advocate for our scheme.

1. Graph on $|V|$ vertices can be used to encode number of $l = \frac{(|V|-1)|V|}{2}$ digits, while $col(D)$ has $|V|$ elements. Hence, overhead information needed to detect errors decreases rapidly with the size of the number to be checked.

2. Proposed scheme is specially suitable for situations, when number and type of errors cannot be easily foreseen. This is major difference between graph coloring based check digit scheme and the others, that successfully detect only certain specified types of error (e.g. single digit error, transpositions, etc., see [5]). In the proposed scheme, detection capabilities and feasibility increase with the volume of information transferred.

3. Proposed scheme can detect all $D'$, such that $col(D)$ is different from $col(D')$. As outlined in the section 2 even single edge change in the graph structure (corresponding to single bit change in the number) can result in change of the graph coloring. When changes in $D$ take place, $D'$ can correspond to any graph from the set of the graphs on $|V|$ vertices (the worst case of error type, that does not allow easy classification and analysis). There are $\Gamma(|V|) = 2^{\frac{|V|(|V|-1)}{2}}$ possible graphs $G(D')$ for the transmitted number $D'$. Let $p$ be the probability that $D \neq D'$ gets undetected through check-digit scheme. Define the probability $p_n$ that $col(D)$ is also proper coloring for $G(D')$, under the condition $D \neq D'$,

$$p_n = \frac{\Gamma(|V|, n)}{\Gamma(|V|)}. \quad (4)$$

The number of $n$-colorable graphs on $|V| = n+y$ vertices decreases very quickly as $y$ increases; from (2) we obtain that $p_n \leqslant 2^{-y}$.

In addition, coloring using $(n-1)$ colors is checked for $D'$. This eliminates all $D \neq D'$ such that $\chi(G(D')) < n$. As the result $p < p_n \leqslant 2^{-y}$ yields $p < 2^{-y}$. Conclusions from upper bound for $\Gamma(|V|, n)$ theorem can be used to derive better estimates for $p$ value.

4. In order to ensure required reliability for the check digit scheme, proper value of $k$ such that $y \geqslant k$, should be maintained for every $D$. This can be achieved by **A** matrix extension. There are numerous techniques that can be used for this purpose. We use example to illustrate the process and its purpose.

---

[1] When $col(D') \neq col(D)$, error would be detected with high probability (same reasoning as for the rest of the scheme applies). For the sake of simplicity we assume that check digits ($col(D)$) are not modified during transmission. For instance, $col(D)$ can be transmitted separately via reliable chanel. Analysis of possible errors in the check digits is not complicated, but this topic is beyond scope of this paper.

Extension is performed by adding to **A**'s bottom additional $k$ rows. These rows have carefully chosen entries and length in such way that resulting matrix remains square. Added entries correspond to extra $k$ vertices added to the graph $G$, resulting in new expanded graph $G_e$. Algorithm for choosing the entries should guarantee that although number of vertices in $G_e$ is $|V| + k$, $\chi(G_e) \leqslant |V|$ (e.g. by using Brook's theorem). Only $D$ and coloring of $G_e$ need to be transmitted. Remaining part of $G_e$ would be recovered using $D'$(algorithm for entries addition is known to both communicating parties). Coloring would be checked for the restored $G_e$.

5. Now is the time to address problem of finding minimal graph coloring. As discussed in the section 2, this problem is NP complete in general case. However, in special cases there are methods for finding it's solution in polynomial time. Hence, when numbers that are to be checked result in graphs from special classes, such methods can be used. Otherwise efficient approximation algorithms, as mentioned in the section 2, have to be used. Instead of looking for minimal graph coloring, approximate solution can be applied. Only consequence of such approach will be reduced value of $y$, compared to minimal coloring case. This is an optimization problem (with: efficiency of coloring approximation algorithm, graph extension algorithm and $y$ as the parameters), that should be solved for particular implementations.

## 5 Concluding remarks and further research

1. In the paper we considered verification method for the bit stream. To make the model more realistic one should join the proposed check-digit scheme with some communication model, which simulates the occurrence of bit errors during transmission process. Additionally, more sophisticated conversion methods between graphs and binary numbers can be considered in order to optimize proposed check-digit scheme.

2. The next step of the research should be an application of the graph coloring check digit scheme to Verifiable Secret Sharing. Check-digit scheme used as part of shares verification procedure, allows to design new approach in this field. Interesting features of such of such approach are:

  a. it is independent on the secret sharing scheme,

  b. it does not matter whether secret shares are numbers or graphs.

3. Another application can be connected with zero knowledge proofs (for instance see[6] ).